\documentclass[final,3p,times]{elsarticle}

\usepackage{algorithm}
\usepackage{algpseudocode}
\usepackage{amsmath}
\usepackage{multirow}
\usepackage{xcolor}
\usepackage{hyperref}
\usepackage{url}

\usepackage{booktabs,threeparttable,caption2}
\usepackage{float}
\usepackage{amsthm,amsmath}
\usepackage{latexsym}
\usepackage{mathrsfs}
\usepackage{mathptmx}
\usepackage{tabularx}
\usepackage{indentfirst}
\biboptions{numbers,sort&compress}

\newdefinition{example}{Experiment}

\newdefinition{rmk}{Remark}

\newcaptionstyle{left}{
\usecaptionmargin\captionfont
{\flushleft\bfseries\captionlabelfont\captionlabel\par}
\mbox{\onelinecaption{\captiontext}{\captiontext}}
}

\begin{document}

\begin{frontmatter}


\title{ Efficient and accurate computation to the $\varphi$-function and its action on a vector} 

\author[author1]{Siyu Yang}
\author[author1,author2]{Dongping Li\corref{cor1}}
\ead{lidp@ccsfu.edu.cn}
\cortext[cor1]{Corresponding author.}
\address[author1]{Department of Mathematics, Changchun Normal University, Changchun 130032, PR China}
\address[author2]{Department of Mathematics, Jilin University, Changchun 130012, PR China}

\begin{abstract}
In this paper, we develop efficient and accurate algorithms for evaluating  $\varphi(A)$ and $\varphi(A)b$, where $A$ is an $N\times N$ matrix, $b$ is an $N$ dimensional vector and $\varphi$ is the function defined by $\varphi(z)\equiv\sum\limits^{\infty}_{k=0}\frac{z^k}{(1+k)!}$.
Such matrix function (the so-called $\varphi$-function) plays a key role in a class of numerical methods well-known as exponential integrators.
The algorithms use the scaling and modified squaring procedure combined with truncated Taylor series. The backward error analysis is presented to find the optimal value of the scaling and the degree of the Taylor approximation. Some useful techniques are employed for reducing the computational cost. Numerical comparisons with state-of-the-art algorithms show that the algorithms perform well in both accuracy and efficiency.
\end{abstract}


\end{frontmatter}

\section{Introduction}
In this work, we consider numerical methods for approximating the first matrix exponential related function and its action on a vector, that is,
\begin{equation}\label{1.1}
\varphi(A)~~\text{and}~~\varphi(A)b,
\end{equation}
 where
\begin{equation}\label{1.2}
\varphi(z)=\sum\limits^{\infty}_{k=0}\frac{z^k}{(1+k)!},  ~A\in \mathbb{C}^{N\times N},~ b\in \mathbb{C}^{N}.
\end{equation}
The $\varphi$-function  satisfies the recursive relation
\begin{equation}\label{1.3}
\varphi(z)=\frac{e^z-1}{z}.
\end{equation}
The problem of numerically approximating such matrix function
is of great importance and is commonly encountered in the solution of constant inhomogeneous
linear system of ordinary differential equations and in the exponential integrators for solving semi-linear problems.
For example, the well-known exponential Euler method for solving the
autonomous semi-linear problems of the form
\begin{equation}\label{1.4}
y'(t)=Ay(t)+N(y(t)),~~y(t_n)=y_n
\end{equation}
yields
\begin{equation}\label{1.5}
y_{n+1} = e^{hA}y_{n-1} + h\varphi(hA)N(y_n).
\end{equation}
If Eq. (\ref{1.4}) has a constant inhomogeneous term, i.e., $N(y(t))\equiv b$, then the scheme (\ref{1.3}) is the exact solution of (\ref{1.4}).
Utilizing the relationship (\ref{1.3}), it can be shown that (\ref{1.5}) is equivalent to
\begin{equation}\label{1.6}
y_{n+1}= h\varphi(hA)(AN(y_n)+y_{n-1}).
\end{equation}
The main cost in the scheme (\ref{1.6}) originates from the need to accurately solve the $\varphi$-function at each time step.
For a detailed overview on exponential integrators, see \cite{Hochbruck2010,BV2005}.

Over the past few years, there has been a tremendous effort to develop efficient
approaches to deal with such matrix functions, see, e.g., \cite{AlMohy2011,Beylkin,Caliari,Hochbruck1998,Sidje1998,Lu,AK,Niesen2012,Skaflestad}.
These methods are generally divided into two classes. The first class of methods compute $\varphi(A)$ explicitly. Among them, the
scaling and modified squaring method combined with Pad\'{e} approximation \cite{Hochbruck1998,Skaflestad} is perhaps the most popular choice for small and medium sized $A$. The method is a variant of the well-known scaling and squaring approach for computing the matrix exponential \cite{AlMohy2009,Higham2005}.
An alternative computation is based on the formula \cite{Saad}:
\begin{equation}\label{1.7}
\begin{aligned}
e^{\mathbb{A}}=\left(\begin{tabular}{cccccc}
$e^{A}$ & $\varphi(A)$\\
$0$ & $I$
\end{tabular}%
\right),
\end{aligned}
~\text {where}~
\begin{aligned}
\mathbb{A}=\left(\begin{tabular}{cccccc}
$A$ & $I$\\
$0$ & $0$
\end{tabular}%
\right)\in \mathbb{C}^{2N\times 2N}.
\end{aligned}
\end{equation}
Thus the computation of $\varphi(A)$ can be reduced to that of the matrix exponential. The effective evaluation of matrix exponential, which arise in many areas of science and engineering, have been extensively investigated in the literature; see, e.g., \cite{AlMohy2009,Defez2018,DP00,Higham2005,Higham,Moler2003,Saad,Sastre2015,Ward} and the references given therein.

In some applications, it requires the computation of matrix-function vector product $\varphi(A)b$ rather than the single $\varphi(A)$.
When $A$ is very large, it is prohibitive to explicitly compute $\varphi(A)$ and then
form the the product with vector $b$. The second class of methods enable evaluation of $\varphi(A)b$ using matrix-vector products and avoids 
the explicit computation of the generally dense matrix $\varphi(A)$.
This type of methods is especially well-suited to large and sparse $A$.
We mention two typical strategies in such an approach: Krylov subspace methods \cite{Sidje1998,Niesen2012} and the scaling-and-squaring method \cite{AlMohy2011}. The former are iterative and difficult to determine a reasonable convergence
criterion to guarantee a sufficiently accurate approximation.
The latter evaluate $\varphi(A)b$ by computing the action of a matrix exponential $e^{\mathbb{A}}$ of dimension $N+1$ on a vector. The method is numerical stable and can achieve a machine accuracy in exact arithmetic.

In the present paper we focus on the direct approach and develop the scaling and modified squaring method in combination with Taylor series to efficiently and accurately evaluate $\varphi(A)$ and $\varphi(A)b,$ respectively. The backward error are used to determine the scaling value $s$ and the Taylor degree $m$. Numerical experiments with other state-of-the-art MATLAB routines illustrate that a straight implementation of the scaling and modified squaring algorithm may be the most efficient.

This paper is organized as follows. Section 2 presents two algorithms for computing $\varphi(A)$. Section 3 deals with algorithm
for evaluating $\varphi(A)b.$ Numerical experiments are given to illustrate the benefits of the algorithms in Section 4. Finally, conclusions are given in Section 5.

Throughout the paper, we use $\|\cdot\|$ to denote an induced matrix norm, and in particular $\|\cdot\|_1$, the 1-norm.
Let $I$ be the identity and 0 be the zero matrix or vector whose dimension are
clear from the context. $e_i$ denotes the $i$-th coordinate vector with appropriate size. $\lfloor x\rfloor$ denotes the largest integer not exceeding $x$ and $\lceil x\rceil$ denotes the smallest integer not less than $x$. Standard MATLAB notations are used whenever necessary.

\section{Computing $\varphi(A)$}
For a given matrix $A\in \mathbb{C}^{N\times N},$ the scaling and modified squaring method exploits the identity \cite{Hochbruck1998,Skaflestad}
\begin{equation}\label{2.1}
\varphi(A)=\frac{1}{2} \varphi(\frac{1}{2} A)(e^{\frac{1}{2}A}+I).
\end{equation}
Applying recursively (\ref{2.1}) $s$ times yields
\begin{eqnarray}\label{2.2}
\varphi(A)=(\frac{1}{2})^{s}\varphi(X)(e^{X}+I)(e^{2X}+I)\ldots(e^{2^{s-1}X}+I),~~s\in \mathbb{N},
\end{eqnarray}
where $X=2^{-s}A$.
The $\varphi(A)$ then can be evaluated using rational polynomial to approximate $\varphi(X)$ and $e^X$ and employing the following
coupled recurrences:
\begin{equation}\label{2.3}
\left\{
\begin{array}{l}
\varphi(2X)= \frac{1}{2} \varphi(X)(e^{X}+I),\vspace{1ex}\\
e^{2X}=e^X\cdot e^X.
\end{array}
\right.
\end{equation}
The scaling parameter $s$ is chosen such that $\|X\|$ is sufficiently small and the method can achieve a prescribed accuracy.

In our algorithm, we use the truncated Taylor series $T_{m}(X)$ to approximate $\varphi(X)$, i.e.,
\begin{equation}\label{2.4}
T_{m}(X) :=\sum\limits^m_{k=0}\frac{X^k}{(1+k)!}.
\end{equation}
Then, the approximation to $e^X$ is naturally chosen as
\begin{equation}\label{2.5}
\tilde{T}_{m}:=XT_{m}(X)+I.
\end{equation}
Here $e^X=\tilde{T}_{m}(X)+\mathcal{O}(X^{m+2})$. The computation of $\tilde{T}_{m}$ only requires one matrix multiplication and one matrix
summation.

In practical the truncated Taylor series $T_m(X)$ in (\ref{2.4}) can be computed efficiently by using the Paterson-Stockmeyer method \cite{Paterson}. The expression is
\begin{equation}\label{2.6}
\tilde{T}_{m}=\sum\limits^r_0B_k\cdot (X^q)^k,~~~r=\lfloor m/q\rfloor,
\end{equation}
where $q$ is a positive integer and
\begin{equation}\label{2.7}
B_k=\left\{
\begin{array}{l}
\sum\limits^{q-1}_{i=0}\frac{1}{(1+qk+i)!} X^i,~~~k=0,1,\ldots,r-1,\vspace{1ex}\\
\sum\limits^{m-qr}_{i=0}\frac{1}{(1+qr+i)!} X^i,~~~k=r.\\
\end{array}
\right.
\end{equation}
Applying Horner's method to (\ref{2.6}), then the number of matrix multiplications for computing $T_m(X)$ is minimized by $q$ either $\lfloor\sqrt {m}\rfloor $ or $\lceil\sqrt {m}\rceil$, and both choices yield the same computational cost.

To obtain a more accurate approximation to $\varphi(X)$, we compute $T_m$ only when $m$ belongs to the index sequence $\mathbb{M}=\{2, 4, 6, 9, 12, 16, 20, 25, 30, 36,\ldots\}$. Assume that $m_i$ is the $i$-th element of the set $\mathbb{M}$, it is shown in \cite[Table 4.1]{Higham}, \cite{Sastre} that the number of matrix multiplications for computing $T_{m_i}(X)$ is the same amount as $T_k(X)$ of $m_{i-1}<k<m_i$.
Then the number of matrix multiplications for computing $T_m(X)$ is
\begin{equation}\label{2.8}
\pi_m=\lceil\sqrt {m}\rceil+\lfloor m/\lceil\sqrt {m}\rceil \rfloor-2.
\end{equation}

Table \ref{tab2.1} lists the corresponding number of matrix multiplications $\pi_m$ to evaluate $T_m$ for the first 12 values of $m$ belonging to $\mathbb{M}$. A brief sketch of the algorithm for solving $\varphi(A)$ is given in Algorithm \ref{alg1}.
\begin{table}[h]
\setlength{\abovecaptionskip}{0.cm}
\setlength{\belowcaptionskip}{-0.3cm}
\caption{Number of matrix multiplications $\pi_m$ required to evaluate $T_m$ for the first 12 optimal values of $m$. }\label{tab2.1}
\begin{center}
\begin{tabular*}{\textwidth}{@{\extracolsep{\fill}}@{~}c|cccccccccccc}
\toprule%
$m$ & $2$ & $4$ & $6$ & $9$ & $12$ & $16$ & $20$ & $25$ & $30$ & $36$ & $42$ & $49$ \\
\hline
$\pi_m$ & $1$ & $2$ & $3$ & $4$ & $5$ & $6$ & $7$ & $8$ & $9$ & $10$ & $11$ & $12$\\
\bottomrule
\end{tabular*}
\end{center}
\end{table}

\begin{algorithm}[htb]
\caption{~Given $A \in \mathbb{C}^{N\times N},$ this algorithm computes $\varphi(A)$ by the scaling and modified squaring based on Taylor series.}\label{alg1}
\begin{algorithmic}[1]
\State Select optimal values of $m$ and $s.$
\State Compute $X=2^{-s}A.$
\State Compute $T=\sum\limits^m_{k=0}\frac{X^k}{(1+k)!}$ by PS method.
\State Compute $\tilde{T}:=XT+I.$
\State Compute $Y=\tilde{T}+I.$
\For{$i=1:s-1$}
\State Compute $\tilde{T}=\tilde{T}^2.$
\State Compute $Y=\frac{1}{2}Y(\tilde{T}+I).$
\EndFor
\State Compute $Y=\frac{1}{2}TY.$
\Ensure~$Y$
\end{algorithmic}
\end{algorithm}

Now we consider the concrete choice of $m$ and $s$. We formulate two approaches to choose the scaling value $s$ and  the Taylor degree $m$, which were similarly introduced in \cite{AlMohy2009,AlMohy2011}. Define the function $h_{m+2}(X)=\log(e^{-X}\tilde{T}_{m}(X))$, then
\begin{equation}\label{2.9}
\tilde{T}_{m}(X)=e^{X+h_{m+2}(X)}
\end{equation}
and
\begin{eqnarray}\label{2.10}
\begin{aligned}
\varphi(A)&\approx(\frac{1}{2})^{s}T_m(\tilde{T}_{m}+I)(\tilde{T}_{m}^{2}+I)\cdots(\tilde{T}_{m}^{2^{s-1}}+I)\\
&=(\frac{1}{2})^{s}\frac{\tilde{T}_{m}-I}{X}(\tilde{T}_{m}+I)(\tilde{T}_{m}^{2}+I)\cdots(\tilde{T}_{m}^{2^{s-1}}+I)\\
&=\frac{e^{2^sX+2^sh_{m+2}(X)}-I}{2^sX}\\
&=\frac{e^{A+\Delta A}-I}{A},
\end{aligned}
\end{eqnarray}
where $\Delta A=2^sh_{m+2}(X)$ is the backward error resulting from the approximation of $\varphi(A).$

Let $X\in \Omega_m:=\{ X\in C^{N\times N}:~~\rho(e^{-X}\tilde{T}_{m}-I)<1\},$ then the function $h_{m+2}(X)$ has a power series expansion
\begin{equation}\label{2.11}
h_{m+2}(X)=\sum\limits^\infty_{k=m+2}c_kX^k.
\end{equation}
By Theorem 4.2 of \cite{AlMohy2009} we have
\begin{equation}\label{2.12}
\frac{\|\Delta A\|}{\|A\|}=\frac{\|h_{m+2}(X)\|}{\|X\|}\leq \tilde{h}_{m+2}(2^{-s}\alpha _p(A)),~~ p(p-1)\leq m+2,
\end{equation}
where $\tilde{h}_{m+2}(x)=\sum\limits^\infty_{k=m+2}|c_k|x^{k-1}$ and $\alpha _p(A)=\max(\|A^p\|^{1/p},\|A^{p+1}\|^{1/(p+1)}).$
Given a tolerance $\text{Tol}$, one can computes
\begin{equation}\label{2.13}
\theta_{m}=max\{\theta:{\tilde{h}_{m+2}(\theta)}\leq \text{Tol}\}.
\end{equation}
Table \ref{tab3.1} presents the maximal values $\theta_m$ satisfying the backward error bound (\ref{2.13}) of $\text{Tol}=2^{-53}$ for the first 12 values of $m$ in $\mathbb{M}$. Thus, once the scaling $s$ is chosen such that
\begin{equation}\label{2.14}
2^{-s}\alpha_p(A) \leq \theta_{m},~~p(p-1)\leq m+2,
\end{equation}
it follows that
\begin{equation}\label{2.15}
\|\Delta A\|\leq\|A\| \cdot\text{Tol}.
\end{equation}
A straightforward computation of inequality (\ref{2.14}) yields
\begin{equation}\label{2.16}
s\geq \lceil \log_2(\alpha _p(A)/\theta_{m}) \rceil.
\end{equation}
We naturally choose the smallest $s$ so that the inequality (\ref{2.14}) holds.
 The total number of matrix multiplications $C_m$ to evaluate $\varphi(A)$ then is
\begin{equation}\label{2.17}
C_m=\pi_m+2s=\lceil \sqrt{m}\rceil+\lfloor \sqrt{m}\rfloor-2+2\max(\lceil \log_2(\alpha _p(A)/\theta_{m}) \rceil,~0 ).
\end{equation}

\begin{table}[h]
\setlength{\abovecaptionskip}{0.cm}
\setlength{\belowcaptionskip}{-0.3cm}
\caption{Maximal values $\theta_m$ such that the backward error bound (\ref{2.13}) does not exceed $\text{tol}=2^{-53}$ for the first 12 optimal values of $m$.}\label{tab3.1}
\begin{center}
\begin{tabular*}{\textwidth}{@{\extracolsep{\fill}}@{~}c|cccccccccccc}
\toprule%
$m$ & $2$ & $4$ & $6$ & $9$ & $12$ & $16$ & $20$ & $25$ & $30$ & $36$ & $42$ & $49$\\
  \hline
  $\theta_m$ & $1.39\text{e-}5$ & $2.40\text{e-}3$ & $2.38\text{e-}2$ & $1.44\text{e-}1$ & $4.00\text{e-}1$ & $9.31\text{e-}1$ & $1.62$ & $2.64$ &$3.77$ &$5.22$&$6.73$&$8.55$ \\
\bottomrule
\end{tabular*}
\end{center}
\end{table}

\begin{figure}[h]
\centering
\includegraphics[width=10cm,height=6cm]{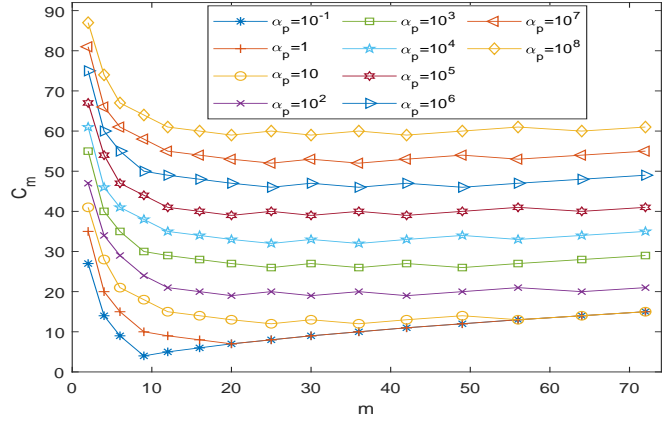}\\
\caption{$m$ versus cost $C_m$ with different $\alpha_p.$}\label{fig2.1}
\end{figure}

In Figure \ref{fig2.1} we have plotted $C_m$
as a function of $m$ for ten different values of $\alpha _p(A)$. We see the location of the first optimal value of $m$, that is, the first value that minimizes $C_m$ is no more than 25. Thus we consider $m$ with $m\in \{2, 4, 6, 9, 12, 16, 20, 25\}$ in the remainder of the section.

In order to get the optimal value of $m$, we consider the following two strategies:

$\bullet$ Choose the first $m\in \{2, 4, 6, 9, 12, 16, 20, 25\}$ such that $\eta_m\leq \theta_m$, where $\eta_m=\min\{\alpha_p(A), p(p-1)\leq m+2\}$, and set $s=0$. When $\eta_{25}>\theta_{25}$, set $m=25$ and $s=\lceil\log_2(\eta_{25}/\theta_{25})\rceil$. To reduce the computational cost, in practical implementation the bound $\|A^p\|^{1/p}$ are estimated using the products of bounds or norms of matrices that have been computed. The details of the process are summarized in Algorithm \ref{alg2}.

$\bullet$ Select the parameters $m$ and $s$ such that the total computational cost (\ref{2.16}) is the lowest. This requires pre-evaluating the first six
1-norm of matrix power, i.e., $\|A^k\|$, $k=1,2,\cdots,6$. The full procedure is given in Algorithm \ref{alg3}.

\begin{algorithm}[htb]
\caption{~Given $A \in \mathbb{C}^{N\times N},$ this algorithm computes the parameters $m,~s$ and $A_i=A^i.$}\label{alg2}
\begin{algorithmic}[1]
\State $s=0.$
\State $A_1=A,$ $A_2=A^2,$ $d_1=\|A_1\|_1,$ $d_2=\|A_2\|_1,$ $d_3=d_1d_2.$
\State $\alpha_1=d_1,$ $\alpha_2=\max(d_2^{1/2},d_3^{1/3}),$ $\eta_1=\alpha_2.$
\If{$\eta_1<=\theta_2,$} $m=2,$ \Return. \EndIf
\If {$\eta_1<=\theta_4,$} $m=4,$ \Return.
\EndIf
\State $A_3=A_1 A_2,$ $d_3=\|A_3\|_1,$ $d_4=\min(d_1d_3,d_4),$ $\alpha_2=\max(d_2^{1/2},d_3^{1/3}),$ $\alpha_3=\max(d_3^{1/3},d_4^{1/4}).$
\State $\eta_2=\min(\alpha_2,\alpha_3).$
\If {$\eta_2<=\theta_6,$} $m=6,$ \Return.
\EndIf
\If {$\eta_2<=\theta_9,$} $m=9,$ \Return.
\EndIf
\State $A_4=A_2^2,$ $d_4=\|A_4\|_1,$ $d_5=\min(d_1d_4,d_2d_3),$ $\alpha_3=\max(d_3^{1/3},d_4^{1/4}),$ $\alpha_4=\max(d_4^{1/4},d_5^{1/5}).$
\State $\eta_3=\min(\alpha_2,\alpha_3,\alpha_4).$
\If {$\eta_3<=\theta_{12},$} $m=12,$ \Return.
\EndIf
\If {$\eta_3<=\theta_{16}$} $m=16,$ \Return.
\EndIf
\State $A_5=A_1A_4,$ $d_5=\|A_5\|_1,$ $d_6=\min(d_1d_5,d_2d_4,d_3^2),$ $\alpha_4=\max(d_4^{1/4},d_5^{1/5}),$ $\alpha_5=\max(d_5^{1/5},d_6^{1/6}).$
\State $\eta_4:=\min(\alpha_2,\alpha_3,\alpha_4,\alpha_5).$
\If {$\eta_4<=\theta_{20}$} $m=20,$ \Return.
\EndIf
\If {$\eta_4<=\theta_{25}$} $m=25,$ \Return.
\EndIf
\State $m=25,$ $s= \lceil\log_2(\eta_4/\theta_{25})\rceil,$
\State $A_i=2^{-is},$ $i=1,\cdots,5.$
\end{algorithmic}
\end{algorithm}

\begin{algorithm}[htb]
\caption{~Given $A \in \mathbb{C}^{N\times N},$ this algorithm computes the parameters $m$ and $s$ based on the number of matrix-matrix products.}\label{alg3}
\begin{algorithmic}[1]
\State $M=[2,~4,~6,~9,~12,~16,~20,~25].$
\State $p_{max}=5,$ $m_{max}=8.$
\State $A_1=A,$ $d_1=\|A\|_1.$
\For {$p = 2:p_{max}+1$}
\State $c =\text {normest}(A,p).$
\State $d_p = c^{1/p}.$
\EndFor
\State $\alpha_1 = d_1.$
\For {$p = 2:p_{max}$}
\State $\alpha_p = \max(d_p,d_{p+1}).$
\EndFor
\State $\eta_1=\alpha_2.$
\For {$p= 2:p_{\max}$}
\State $\eta_p=\min(\eta_{p-1},\alpha_p).$
\EndFor
\For {$m=[2,~4,~6,~9,~12,~16,~20,~25]$}
\If {$m=2$}
\State $s_m=\max(\lceil \log_2(\eta_2/\theta_m)\rceil,0).$
\ElsIf {$m=[4,6,9]$}
\State $s_m=\max(\lceil \log_2(\eta_3/\theta_m)\rceil,0).$
\ElsIf {$m=[12,16]$}
\State $s_m=\max(\lceil \log_2(\eta_4/\theta_m)\rceil,0).$
\Else
\State $s_m=\max(\lceil \log_2(\eta_5/\theta_m)\rceil,0).$
\EndIf
\State $q_m=\sqrt{m},$ $C_m=\lceil q_m\rceil+\lfloor q_m\rfloor-2+2*s_m.$
\EndFor
\State $m=\text{argmin}_{m\in M}C_m,$ $s=s_m,$ $q=\lceil q_m\rceil.$
\State $A_i=A_{i-1}A,$ $i=2:q.$
\State  $A_i=2^{-is}A_i,$ $i=2:q.$
\end{algorithmic}
\end{algorithm}

\section{Computing $\varphi(A)b$}
We now focus our attention on accurately and efficiently evaluating $\varphi(A)b$ for sparse and large matrix $A$.
Following an
idea of Al-Mohy and Higham \cite{AlMohy2011}, we will use the scaling part of the scaling and modified squaring method in combination with truncated Taylor series to approximate the function. The computational cost of the method is dominated by matrix-vector products.

 We start by recalling the following general recurrence \cite{Skaflestad}:
\begin{equation}\label{3.0}
\varphi((\alpha+\beta)z)=\frac{1}{(\alpha+\beta)}[\beta e^{\alpha z}\varphi(\beta z)+\alpha\varphi(\alpha z)],~~~\alpha,\beta\in \mathbb{R},~z\in \mathbb{C}.
\end{equation}
As a special case of (\ref{3.0}), we have
\begin{equation}\label{3.0a}
\varphi(sz)=\frac{1}{s}[e^{(s-1)z}\varphi(z)+(s-1)\varphi((s-1)z)],~~s\in \mathbb{N}.
\end{equation}
Taking $Y=\frac{1}{s}A,$ and using (\ref{3.0a}) it follows that
\begin{eqnarray}\label{3.1}
\begin{aligned}
\varphi(A)=&\frac{1}{s}[e^{(s-1)Y}\varphi(Y)+(s-1)\varphi((s-1)Y)]\\
=&\frac{1}{s}[e^{(s-1)Y}\varphi(Y)+e^{(s-2)Y}\varphi(Y)+(s-2)\varphi((s-2)Y)]\\
&\vdots\\
=&\frac{1}{s}\varphi(Y)[e^{(s-1)Y}+e^{(s-2)Y}+\cdots+e^{2Y}+e^{Y}+I].
\end{aligned}
\end{eqnarray}

Choose the integers $m$ and $s$ such that $\varphi(Y)$ and $e^Y$ can be well-approximated by the truncated Taylor series $T_m(Y)$ and $\tilde{T}_m(Y):=YT_m(Y)+I$ defined by (\ref{2.4}) and (\ref{2.5}).
Then $\varphi(A)b$ can be approximated by firstly evaluating the recurrence
\begin{eqnarray}\label{3.2}
b_1=T_m(Y)b ~~~\text{and}~~~b_{i+1}=\tilde{T}_m(Y)b_i,~~i=1,2,\cdots,s-1,
\end{eqnarray}
and then computing $\frac{1}{s}$ times the sum of $b_i, i=1,2,\cdots,s-1.$
This process requires $s$ multiplications of matrix polynomial with a vector, $s$ vector additions,
and 1 scalar multiplication. The number of matrix-vector products for evaluating $\varphi(A)b$ by recurrence (\ref{3.2}) is
$C_m=s(m+1)-1.$

\begin{algorithm}[htb]
\caption{~Given $A \in \mathbb{C}^{N\times N},~b\in \mathbb{C}^{N\times n_0},$ this algorithm computes $\varphi(A)b$ by the scaling and modified squaring based on Taylor series.}\label{alg4}
\begin{algorithmic}[1]
\State Select optimal values of $m$ and $s.$
\State Compute $Y=A/s.$
\State Compute $b_1=\sum\limits^m_{k=0}\frac{Y^k}{(1+k)!}b$ based on matrix-vector products.
\State Compute $f=b_1.$
\For{$i=1:s-1$}
\State Compute $b_{i+1}=\sum\limits^{m+1}_{k=0}\frac{Y^k}{k!}b_i$ based on matrix-vector products.
\State Compute $f=f+b_i.$
\EndFor
\State Compute $f=\frac{1}{s}f.$
\Ensure~$f$
\end{algorithmic}
\end{algorithm}

The procedure described above mentions two key parameters: the degree $m$ of the matrix polynomial $T_m(Y)$ and the scaling parameter $s$. We use the backward error analysis combined with the computational cost to choose an optimal parameters $m$ and $s$.
The backward error analysis of the method is exactly the same as the above section. The only difference is the form of their scaling coefficients. The former is $2^{-s}$ and the latter is $\frac{1}{s}.$ The relative backward error of the method satisfies
\begin{equation}\label{3.3}
\frac{\|\Delta A\|}{\|A\|}\leq \tilde{h}_{m+2}(\frac{1}{s}\alpha _p(A)),
\end{equation}
where $\tilde{h}_{m+2}(x)$ and $\alpha _p(A)$ are defined exactly as in (\ref{2.8}).
Given a tolerance $\text{Tol}$ and integer $m,$ the parameter $s$ is chosen so that $s^{-1}\alpha _p(A)\leq \theta_m,$ i.e.,
\begin{equation}\label{3.4}
s\geq\lceil\alpha _p(A)/\theta_m\rceil.
\end{equation}
And the cost of the algorithm in matrix-vector products is
\begin{equation}\label{3.5}
C_m=(m+1)\lceil\alpha _p(A)/\theta_m\rceil-1.
\end{equation}

Let $p_{max}$ denote the largest positive integer $p$ such that $p(p-1)\leq m_{max}+2.$ Then the optimal cost is
\begin{equation}\label{3.6}
C_{m^*}=\min\{(m+1)\lceil\alpha _p(A)/\theta_m\rceil-1:~~2\leq p\leq p_{max}, p(p-1)-2\leq m \leq m_{max}\},
\end{equation}
where $m^*$ denotes the smallest value of $m$ at which the minimum is attained. The optimal scaling parameter is $s=C_{m^*}/m^*.$

The cost of computing $\alpha_p(A)$ for $p = 2: p_{max}$ is approximately $2lp_{max}(p_{max} +3),~~l=1~\text{or}~2$. If the cost $C_{m_{max}}$ matrix-vector products of evaluating $\varphi(A)b$ with $m$ determined by using $\|A\|_1$ in place of $\alpha _p(A)$ in (\ref{3.5}) is no larger than the
cost of computing the $\alpha _p(A)$, i.e.
\begin{equation}\label{3.7}
\|A\|_1 <=[4*p_{max}*(p_{max} + 3)+1]*[\theta_{m_{max}}/(m_{max}+1)].
\end{equation}
Then we should certainly use $\|A\|_1$ in place of the $\alpha _p(A)$.
The details of the method is summarized in Algorithms \ref{alg4} and \ref{alg5}.

\begin{algorithm}[htb]
\caption{~Given $A \in \mathbb{C}^{N\times N},~b\in \mathbb{C}^{N\times n_0},$ $p_{max}$ and $m_{max}.$ this algorithm computes the parameters $m$ and $n$ based on the number of matrix-vector products.}\label{alg5}
\begin{algorithmic}[1]
\State $M=1:m_{max}.$
\State $d_1=\|A\|_1.$
\If {$d_1 <=\theta_{m_{max}}*(4*p_{max}*(p_{max}+ 3)+1)/m_{max}$}
\State $m=\arg\min \{(m+1)\lceil d_1/\theta_m\rceil-1,~1\leq m\leq m_{max}\}.$
\State $n=\lceil d_1/\theta_m\rceil,$ \Return.
\EndIf
\For {$p = 2:p_{max}+1$}
\State $c =\text {normest}(A,p).$
\State $d_p = c^{1/p}.$
\EndFor
\State $\alpha_1 = d_1.$
\For {$p = 2:p_{max}$}
\State $\alpha_p = \max(d_p,d_{p+1}).$
\EndFor
\State $[m,p]=\arg\min\{(m+1)\lceil\alpha _p/\theta_m\rceil-1:2\leq p\leq p_{max}, p(p-1)-2\leq m \leq m_{max}\}.$
\State $s=\max(\lceil\alpha _{p}/\theta_{m}\rceil,1).$
\Ensure~$m, s.$
\end{algorithmic}
\end{algorithm}

\section{Numerical experiments}
In this section we perform two numerical experiments to test the performance of the approach that has been presented in the previous sections. All tests are performed under Windows 10 and MATLAB R2018b running on a laptop with an Intel Core i7 processor with 1.8 GHz and RAM 8 GB.

We use Algorithm \ref{alg1} in combination with Algorithm \ref{alg2} and Algorithm \ref{alg3} to evaluate $\varphi(A),$ and Algorithm \ref{alg4} combined with Algorithm \ref{alg5} to compute $\varphi(A)b.$ The three combined algorithms are denoted as \textsf{phitay1}, \textsf{phitay2} and \textsf{phimv}, respectively.

\begin{example} \label{exa1}
In this experiment we compare algorithms \textsf{phitay1} and \textsf{phitay2} with existing MATLAB routine \textsf{phipade13} from \cite{Skaflestad}. The function \textsf{phipade13} employs scaling and modified squaring method based on [13/13] pad\'e approximation to evaluate $\varphi(A).$
We use a total of 201 matrices divided into two sets to test these algorithms. The two test sets are described as follows:

$\bullet$ The first test set contains 62 $8\times 8$ test matrices as in \cite{Higham2005} and \cite[sec. 4.1]{Sastre2015}. The first 48 matrices are obtained from the subroutine \textsf{matrix}\footnote{The subroutine \textsf{matrix} can generate fifty-two matrices. Matrices 17, 42, 44, 43 are excluded the scope of the test as the first three overflow in double precision and the last is repeated as matrix 49.} in the matrix computation toolbox \cite{Highamtool}. The other fourteen test matrices of dimension $2-20$ come from \cite[Ex. 2]{Higham2003}, \cite[Ex. 3.10]{DP00}, \cite[p. 655]{KL1998},  \cite[p. 370]{NH1995}, \cite[Test Cases 1–4]{Ward}.

$\bullet$ The second test set is essentially the same tests as in \cite{Defez2018}, which consists of 139 matrices of dimension $n=128.$ The first 39 matrices are obtained from MATLAB routine \textsf{matrix} of the Matrix Computation Toolbox \cite{Highamtool}. The remaining 100 matrices are generated
randomly, half of which are diagonalizable and half non diagonalizable matrices.

In this implementation, we evaluated the relative errors in the 1-norm of the computed solutions $Y$, i.e.,
\begin{equation}\label{4.1}
Error= \frac{\|Y-\varphi(A)\|_1}{\|\varphi(A)\|_1}.
\end{equation}
The "exact" $\varphi(A)$ is computed using MATLAB build-in function \textsf{expm} of \cite{Higham2005, AlMohy2009} by evaluating the augmented matrix exponential (\ref{1.7}) at 100-digit precision using MATLAB’s Symbolic Math Toolbox.
\end{example}

\begin{figure}[H]
\begin{minipage}{0.5\linewidth}
\centering
\includegraphics[width=7cm,height=5cm]{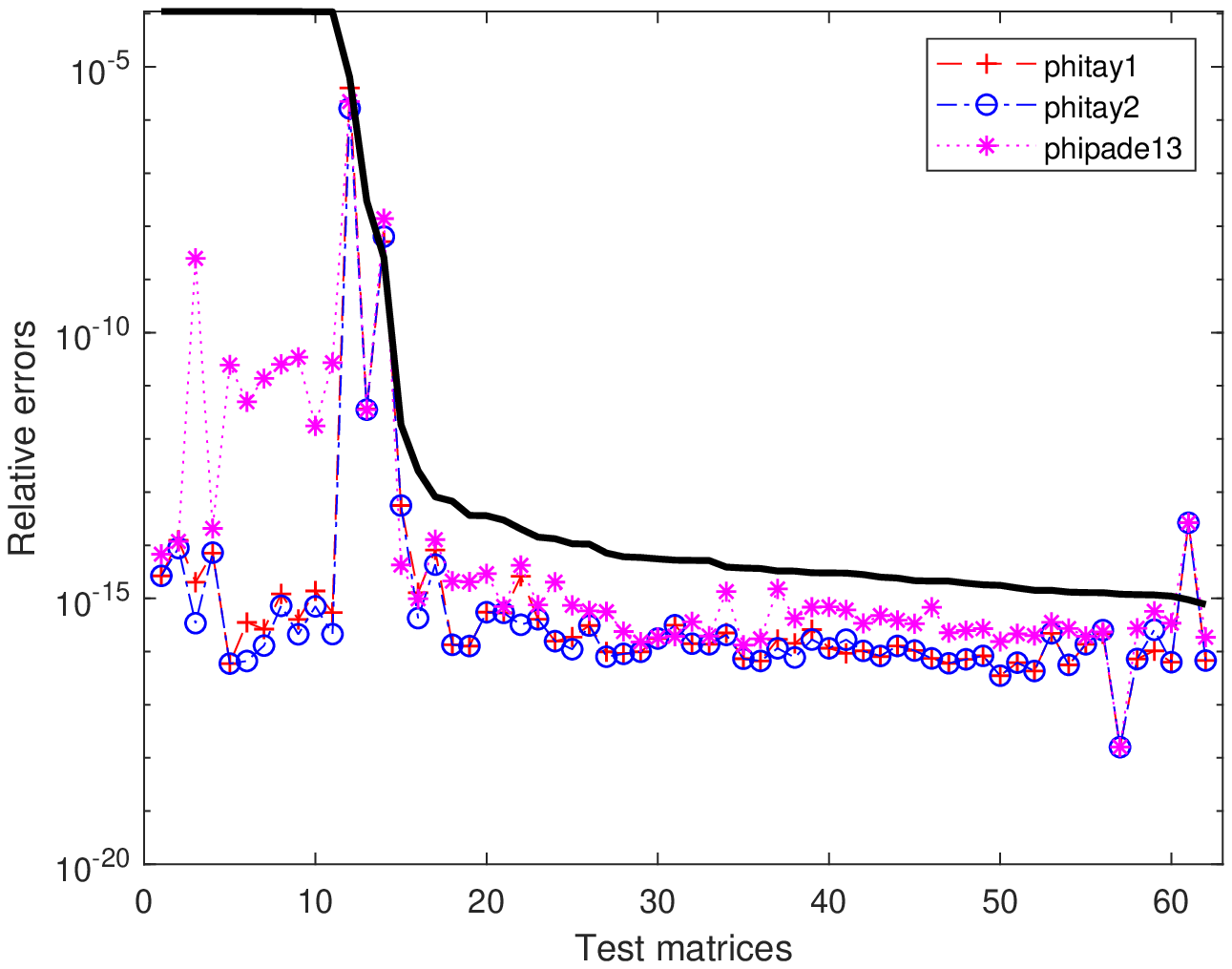}\\
\footnotesize{a. ~Normwise relative errors}
\end{minipage}
\mbox{\hspace{-2cm}}
\begin{minipage}{0.5\linewidth}
\centering
\includegraphics[width=7cm,height=5cm]{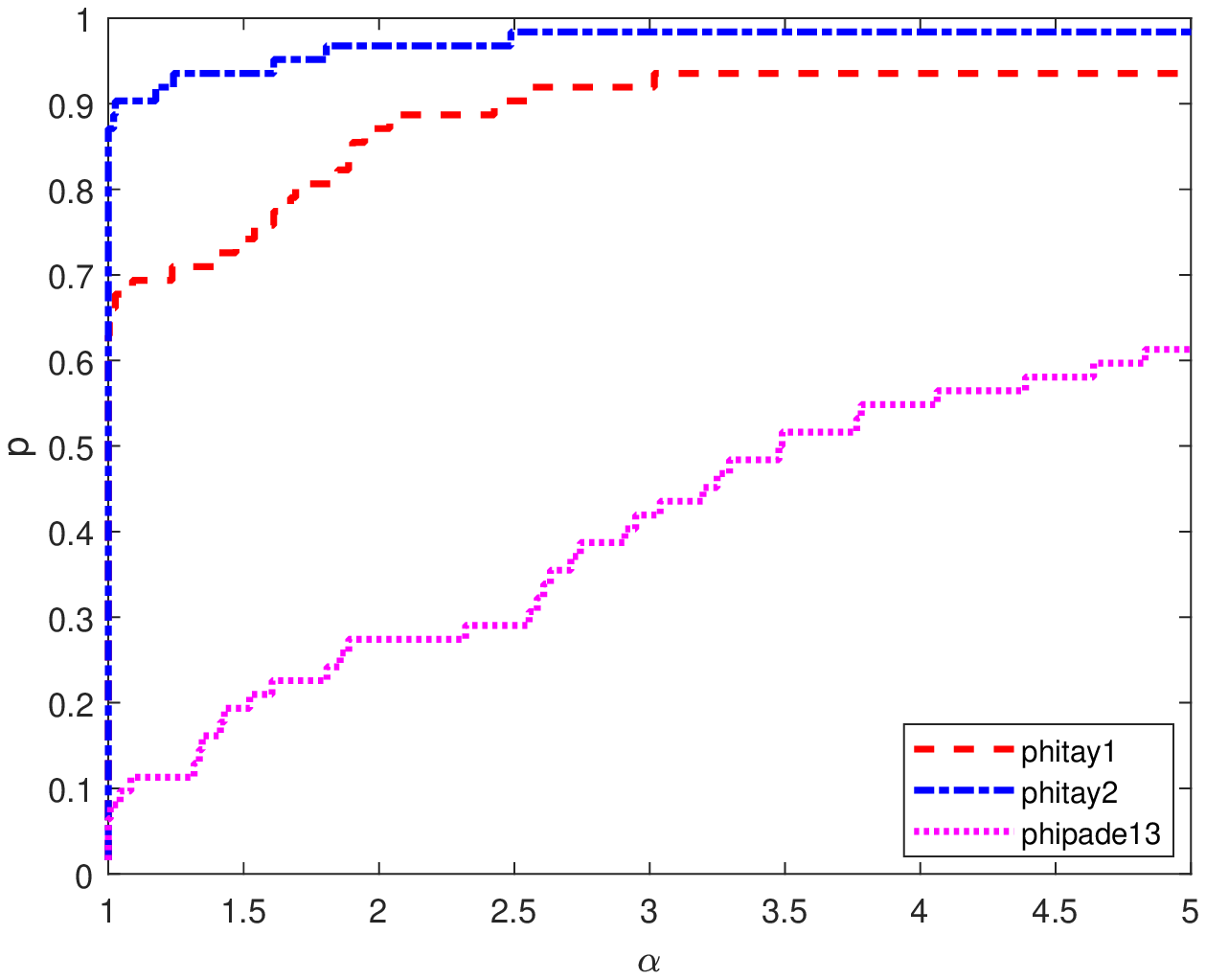}\\
\footnotesize{b.~Performance of errors}
\end{minipage}\\
\begin{minipage}{0.5\linewidth}
\centering
\includegraphics[width=7cm,height=5cm]{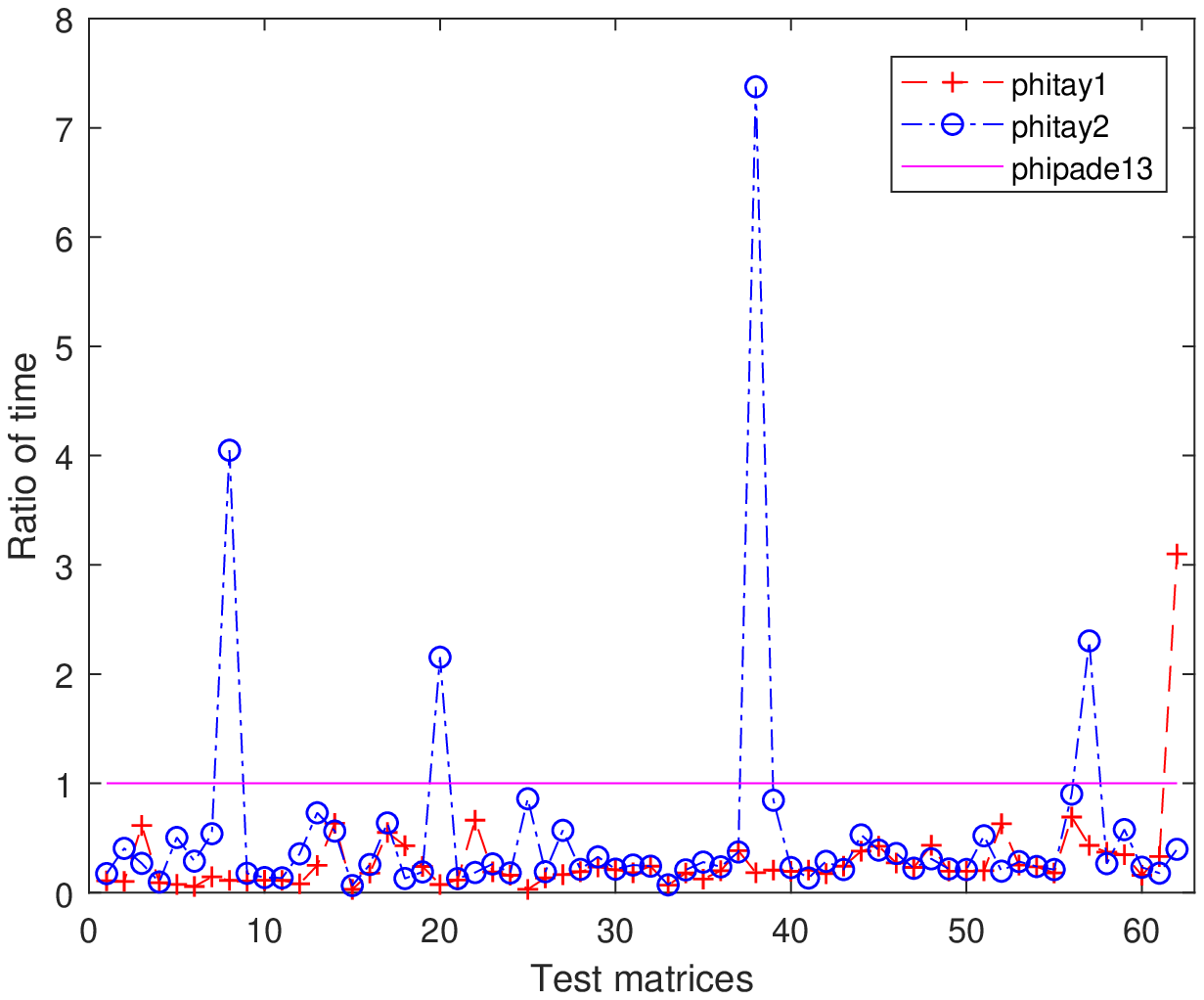}\\
\footnotesize{c.~Ratio of execution times}
\end{minipage}
\mbox{\hspace{-2cm}}
\begin{minipage}{0.5\linewidth}
\centering
\includegraphics[width=7cm,height=5cm]{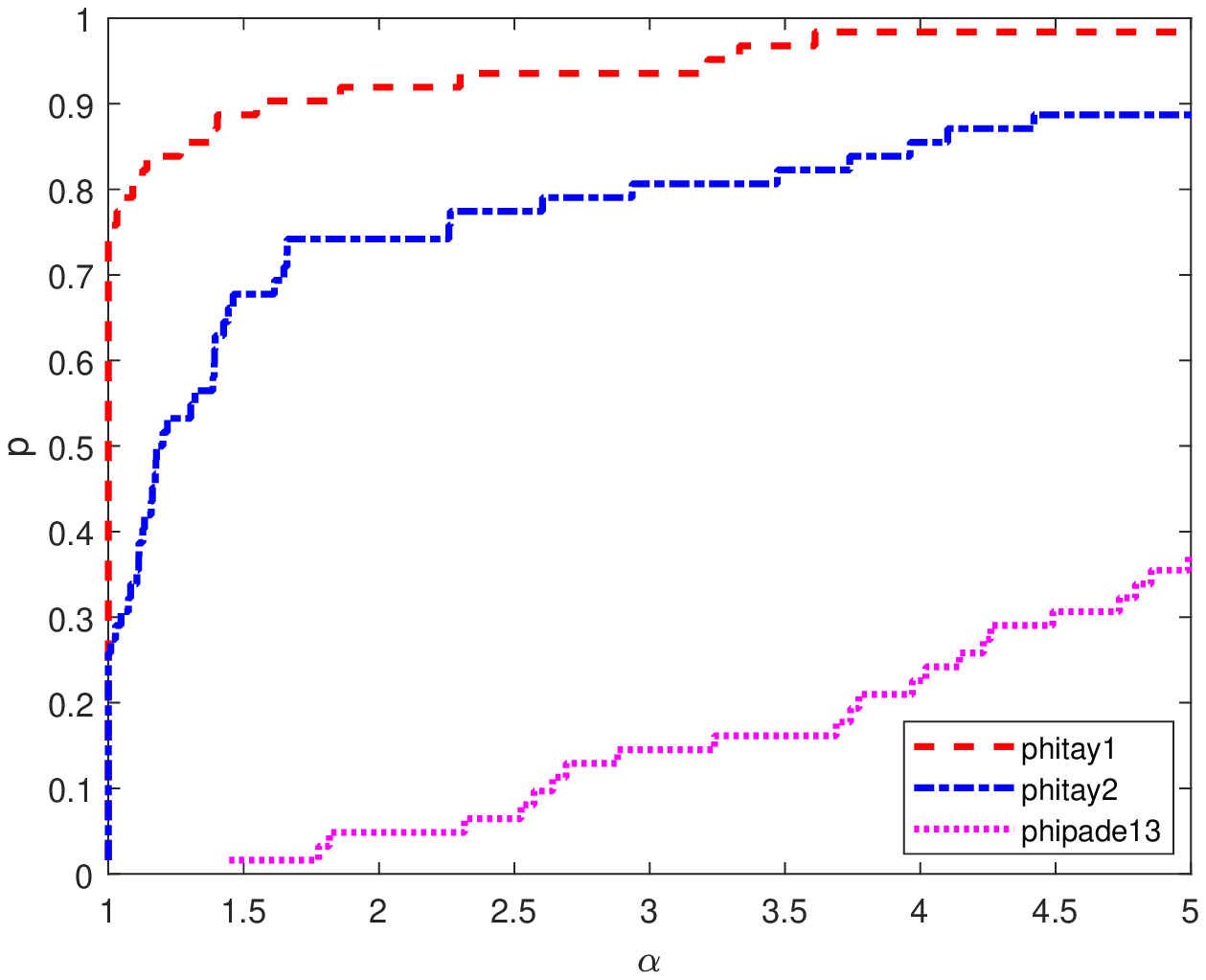}\\
\footnotesize{d.~Performance of execution times}
\end{minipage}
\caption{Results for test matrix set 1 for Experiment \ref{exa1}.}\label{fig4.1}
\end{figure}
\begin{figure}[H]
\begin{minipage}{0.5\linewidth}
\centering
\includegraphics[width=7cm,height=5cm]{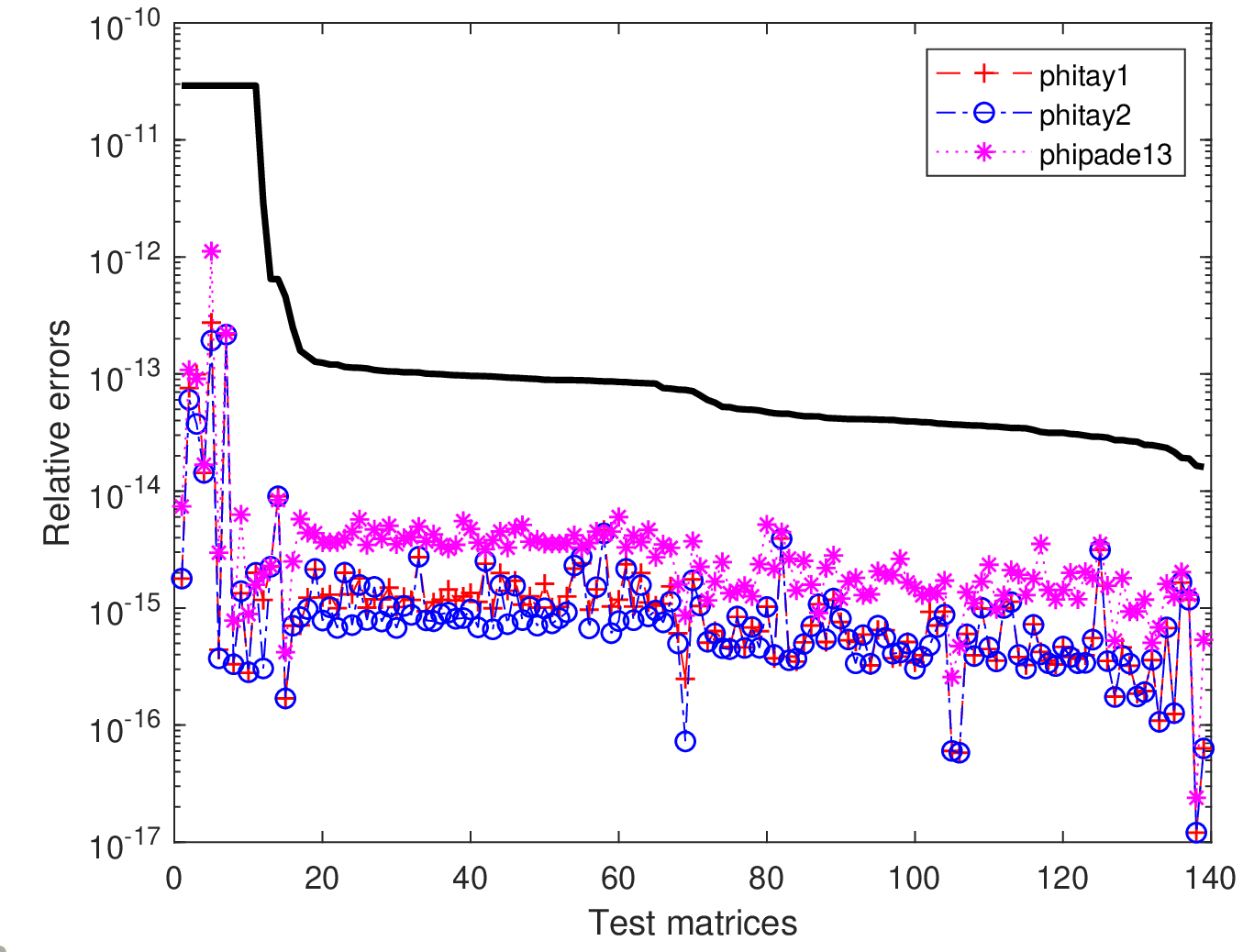}\\
\footnotesize{a. ~Normwise relative errors}
\end{minipage}
\mbox{\hspace{-2cm}}
\begin{minipage}{0.5\linewidth}
\centering
\includegraphics[width=7cm,height=5cm]{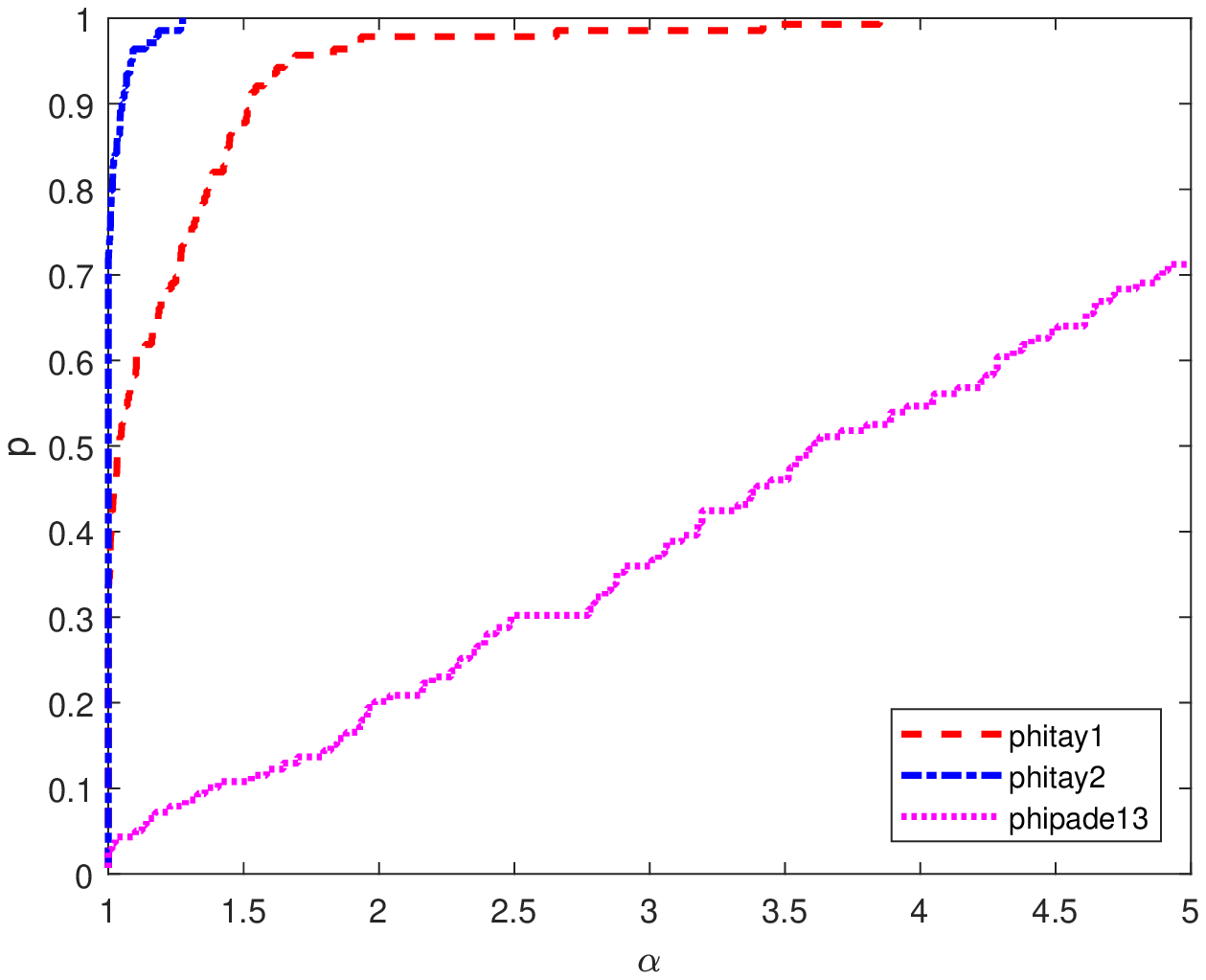}\\
\footnotesize{b.~Performrnce of errors}
\end{minipage}\\
\begin{minipage}{0.5\linewidth}
\centering
\includegraphics[width=7cm,height=5cm]{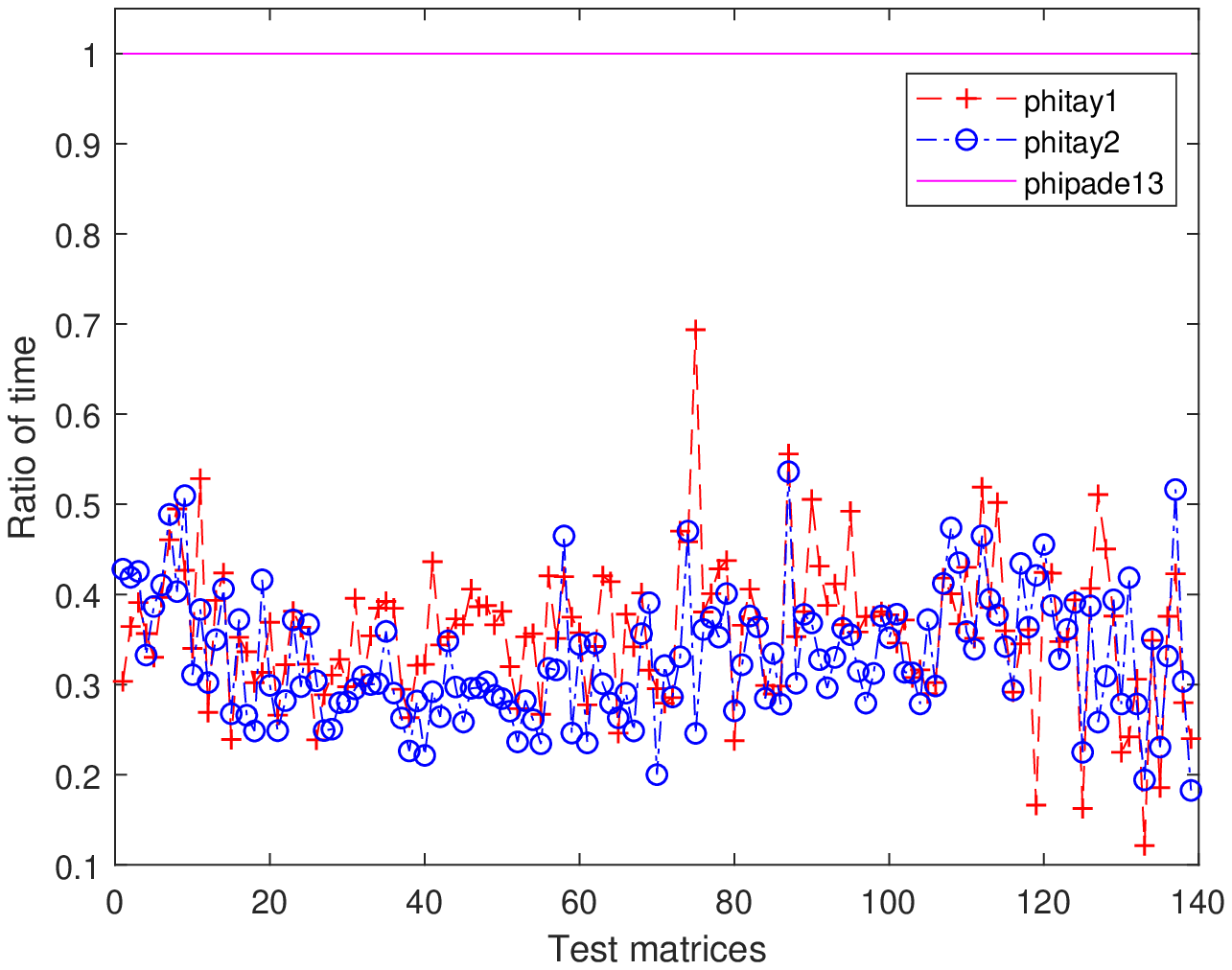}\\
\footnotesize{c.~Ratio of execution times}
\end{minipage}
\mbox{\hspace{-2cm}}
\begin{minipage}{0.5\linewidth}
\centering
\includegraphics[width=7cm,height=5cm]{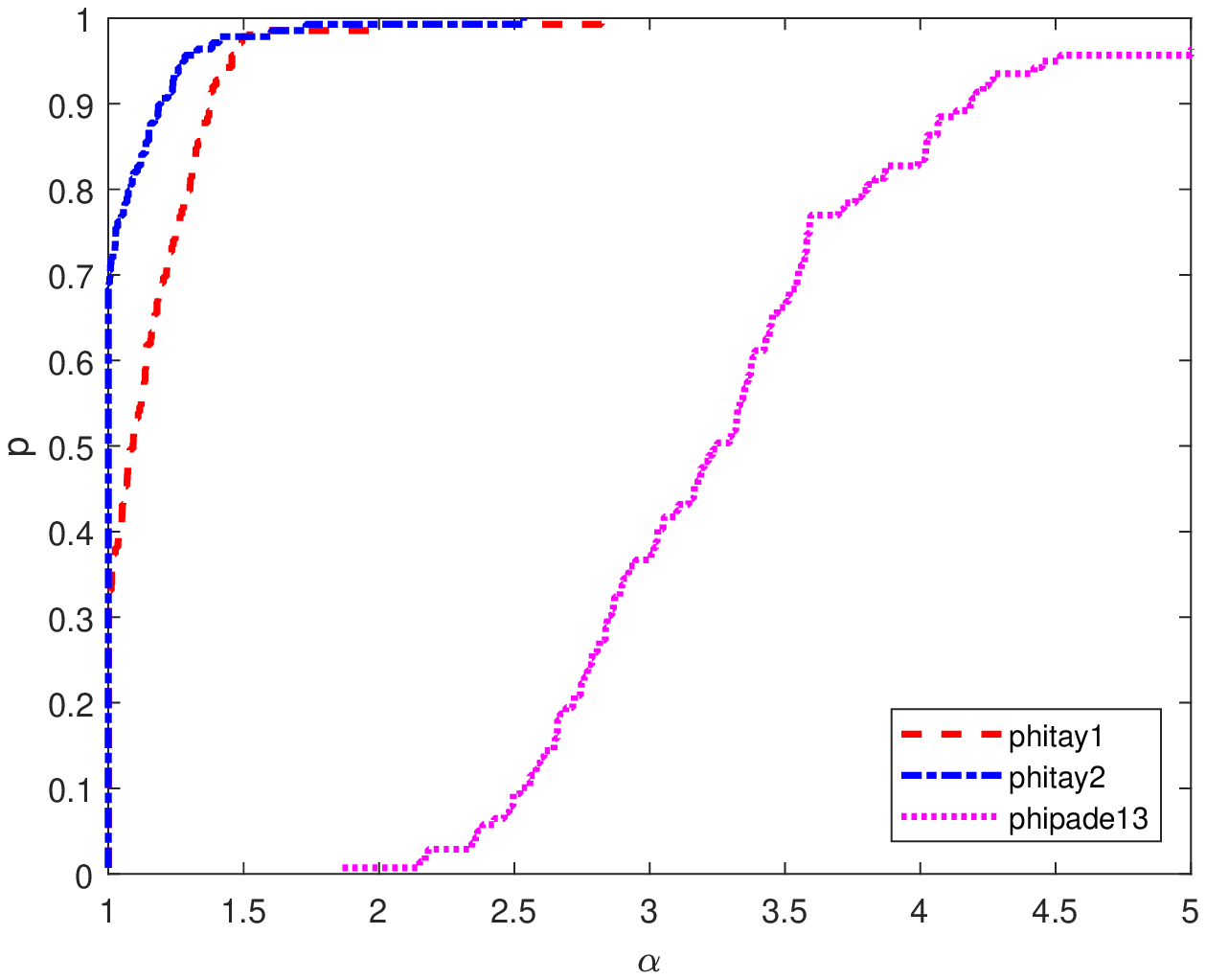}\\
\footnotesize{d.~Performance of execution times}
\end{minipage}
\caption{Experimental results for test matrix set 2 for Experiment \ref{exa1}.}\label{fig4.2}
\end{figure}

In Figs. \ref{fig4.1} and \ref{fig4.2}, we present the relative errors, the performances of relative errors, the ratio of execution times and the performances of execution times for each test. Figs. \ref{fig4.1}(a) and \ref{fig4.2}(a) display the relative error of the algorithms in
our test sets, sorted by decreasing relative condition number of $\varphi(A)$ at $A$. The solid black line represents the unit roundoff multiplied by the relative condition number, which is estimated by the MATLAB routine \textsf{funm\underline{~}condest1} in the Matrix Function Toolbox \cite{Highamtool}. Figs. \ref{fig4.1}(b) and \ref{fig4.2}(b) show the performance profiles of the three solvers on the same error data. For a given $\alpha,$ the corresponding value of $p$ on each performance curve is the probability that the algorithm has a relative error lower than or equal to
$\alpha$ times the smallest error over all the methods involved \cite{DM02}. The results show that the methods based on Taylor series are more accurate than the implementation based
on Pad\'{e} series, and \textsf{phitay2} is slightly more accurate than \textsf{phitay1}. Figs. \ref{fig4.1}(c) and \ref{fig4.2}(c) show the the ratio of execution times of the three solvers with respect to \textsf{phipade13}. The performance on the execution times
of three methods is compared in  Figs. \ref{fig4.1}(d) and \ref{fig4.2}(d). We notice that \textsf{phitay1} and \textsf{phitay2} have lower execution times than \textsf{phipade13}, and the execution time of \textsf{phitay1} is slightly lower than the execution
time of \textsf{phitay2} on test set 1 but the opposite is true on test set 2. This can be attributed to the choice of the key parameters $m$ and $s$ in the methods. The \textsf{phitay2} uses a minimum amount of computational costs to determine the optimal parameters by evaluating exactly the 1-norm of  $\|A^k\|_1^{1/k}$ for a few values of $k$ using a matrix norm estimator. Although this requires some extra calculations in estimating the norm of matrix power, the computational advantages will gain as the dimension of the matrix increases.

\begin{example} \label{exa2}This experiment uses the same tests as \cite{Niesen2012}. There are four different sparse matrices test matrices. The matrices details are

$\bullet$ The first matrix \textsf{orani678} is an unsymmetric sparse matrix of order $N = 2, 529$ with $nnz = 90, 158$ nonzero
elements and its 1-norm is 1.04e+03.

$\bullet$ The second matrix \textsf{bcspwr10} is a symmetric Hermitian sparse matrix of order $ N= 5, 300$ with $nnz = 21,842$ nonzero
elements and its 1-norm is 14.

$\bullet$  The third matrix \textsf{gr\underline{~}30\underline{~}30} is an symmetric sparse matrix of order $N = 900$ with
$nnz = 7, 744$ nonzero elements and its 1-norm is 16.

$\bullet$  The fourth matrix \textsf{helm2d03} is a sparse matrix of order $N = 392, 257$ has $nnz = 2, 741, 935$ nonzero elements and its 1-norm is 10.72.

We use our algorithm \textsf{phimv} with other two popular MATLAB routines \textsf{phiv} of \cite{Sidje1998} and \textsf{phipm} of \cite{Niesen2012} to evaluate $\varphi(tA)b$ and $\varphi_0(tA)b_0+t\varphi_1(tA)b_1$ with $t=10$ for the first test matrix and $t=2$ for the other three, respectively. As in \cite{Niesen2012}, we choose the vectors $b=b_0=b_1=[1, 1, . . . , 1, 1]^T$ except for the second test matrix with $b = [1, 0, . . . , 0, 1]^T.$
The MATLAB routines \textsf{phiv} and \textsf{phipm} are run with their default parameters and the uniform convergence tolerance \textsf{Tol}=\textsf{eps} (\textsf{eps} is the unit roundoff) in our experiments.

In this tests, we assess the accuracy of the computed solution $y$ by the relative errors
\begin{equation}\label{4.2}
Error= \frac{\|y-y_{exa}\|_2}{\|y_{exa}\|_2}.
\end{equation}
 where $y_{exa}$ is a reference solution obtained by computing the action of the augmented matrix exponential using MATLAB routine \textsf{expmv} of \cite{AlMohy2011}. We measure the average ratio of execution times ($t_{ratio}$) of each of the codes relative to \textsf{phikmv} by running the comparisons 100 times. Tables \ref{tab4.1} and \ref{tab4.2} show the numerical results. All three methods deliver almost the same accuracy, but \textsf{phikmv} performs
to be the fastest.
\end{example}

\begin{table}[H]
\setlength{\abovecaptionskip}{0.cm}
\setlength{\belowcaptionskip}{-0.3cm}
\caption{Comparisons of the average speedup of \textsf{phiv}, \textsf{phip} and
\textsf{phipm} with respect to \textsf{phimv} and the relative errors for solving $\varphi(A)b$.}
 \label{tab4.1}
\begin{center}
\begin{tabular*}{\textwidth}{@{\extracolsep{\fill}}@{~~}c|lr|lr|lr|lr}
\toprule%
\raisebox{-2.00ex}[0cm][0cm]{method}&
 \multicolumn{2}{c|}{\textsf{orani678}}&\multicolumn{2}{c|}{\textsf{bcspwr10}}&\multicolumn{2}{c|}{\textsf{gr\underline{~}30\underline{~}30}}&\multicolumn{2}{c}{\textsf{helm2d03}}
\cr
\cmidrule(lr){2-9}
 &error&$t_{ratio}$&error&$t_{ratio}$ &error&$t_{ratio}$ &error&$t_{ratio}$\\
\cmidrule(lr){1-9}
\multirow{1}{*}{\textsf{phiv}}&1.1966e-15&1.38   &4.8688e-16&45.90   &1.8016e-14&9.92   &1.0842e-13&20.98\\
\cmidrule(lr){1-9}
\multirow{1}{*}{\textsf{phipm}}&1.7612e-15&0.79 & 6.6655e-16&7.4  &3.7391e-15&5.54   &1.3550e-13&1.59\\
\cmidrule(lr){1-9}
\multirow{1}{*}{\textsf{phimv}}&1.1682e-15&1  &3.6051e-16&1 &1.2622e-15&1   &6.2692e-14&1\\
\bottomrule
\end{tabular*}
\end{center}
\end{table}

\begin{table}[H]
\setlength{\abovecaptionskip}{0.cm}
\setlength{\belowcaptionskip}{-0.3cm}
\caption{Comparisons of the average speedup of \textsf{phiv}, \textsf{phip} and
\textsf{phipm} with respect to \textsf{phimv} and the relative errors for solving $\varphi_0(tA)b_0+t\varphi_1(tA)b_1$.}
 \label{tab4.2}
\begin{center}
\begin{tabular*}{\textwidth}{@{\extracolsep{\fill}}@{~~}c|lr|lr|lr|lr}
\toprule%
\raisebox{-2.00ex}[0cm][0cm]{method}&
 \multicolumn{2}{c|}{\textsf{orani678}}&\multicolumn{2}{c|}{\textsf{bcspwr10}}&\multicolumn{2}{c|}{\textsf{gr\underline{~}30\underline{~}30}}&\multicolumn{2}{c}{\textsf{helm2d03}}
\cr
\cmidrule(lr){2-9}
 &error&$t_{ratio}$&error&$t_{ratio}$ &error&$t_{ratio}$ &error&$t_{ratio}$\\
\cmidrule(lr){1-9}
\multirow{1}{*}{\textsf{phiv}}&4.7206e-16&1.06   &7.4407e-16&43.63   &2.1790e-15&15.58   &1.5239e-14&20.73\\
\cmidrule(lr){1-9}
\multirow{1}{*}{\textsf{phipm}}&1.0408e-15&0.63 & 1.2216e-15&8.73  &7.4383e-15&5.58   &8.7267e-15&1.77\\
\cmidrule(lr){1-9}
\multirow{1}{*}{\textsf{phimv}}&1.8024e-15&1  &7.6561e-16&1 &8.7257e-16&1   &8.7682e-15&1\\
\bottomrule
\end{tabular*}
\end{center}
\end{table}
\section{Conclusion}
The computation of $\varphi$-functions can lead to a large computational burden of exponential integrators. In this work three accurate algorithms \textsf{phitay1}, \textsf{phitay2} and \textsf{phimv} have been developed to compute the first $\varphi$-function and its action on a vector. The first two are used for solving $\varphi(A)$ and the last one is used for $\varphi(A)b$. These algorithms employ the scaling and modified squaring procedure based on the truncated Taylor series of the $\varphi$-function and are backward stable in exact arithmetic. For \textsf{phitay1} and \textsf{phitay2}, the optimal Horner and Paterson-Stockmeyer’s technique has been applied to reduce the computational cost. The main difference of both is the estimation of matrix powers $\|A^k\|_1^{1/k}$ for a few values of $k$ and \textsf{phitay2} allows to determine the optimal values of scaling and the degree of the Taylor approximation by the minimum amount of computational costs. The \textsf{phimv} takes the similar approach as \textsf{phitay2} to determine the key parameters. The computational costs mostly focused on computing matrix-vector products which is especially well-suited large sparse matrix. Numerical comparisons with other state-of-the-art MATLAB routines illustrate that the methods proposed are efficient and reliable. In the future we hope to further generalize these methods to the general exponential related functions and their linear combination.

\section*{Acknowledgements}
This work was supported in part by the Jilin Scientific and Technological Development Program (Grant Nos. 20200201276JC and 20180101224JC) and the Natural Science Foundation of Jilin Province (Grant No. 20200822KJ), and the Scientific Startup Foundation for Doctors of Changchun Normal University (Grant No. 002006059).

\end{document}